\documentclass[12pt]{amsart}
\input amssym.def
\baselinestretch
\renewcommand{\thefootnote}

\newcommand{\lr}{\mbox{$\longrightarrow$}}

\newcommand{\be}{\begin{equation}}
\newcommand{\ee}{\end{equation}}

\newtheorem{guess}{Theorem}[section]
\newcommand{\bth}{\begin{guess}$\!\!\!$~}
\newcommand{\eeth}{\end{guess}}

\newtheorem{propo}[guess]{Proposition}
\newcommand{\bpropo}{\begin{propo}$\!\!\!$~}
\newcommand{\epropo}{\end{propo}}

\newtheorem{lema}[guess]{Lemma}
\newcommand{\blem}{\begin{lema}$\!\!\!$~}
\newcommand{\elem}{\end{lema}}

\newtheorem{defe}[guess]{Definition}
\newcommand{\bdefe}{\begin{defe}$\!\!\!$~}
\newcommand{\edefe}{\end{defe}}

\newtheorem{coro}[guess]{Corollary}
\newcommand{\bcor}{\begin{coro}$\!\!\!$~}
\newcommand{\ecor}{\end{coro}}

\newtheorem{rema}[guess]{Remark}
\newcommand{\brem}{\begin{rema}$\!\!\!$~\rm}
\newcommand{\erem}{\end{rema}}

\newtheorem{exam}[guess]{Example}
\newcommand{\beg}{\begin{exam}$\!\!\!$~\rm}
\newcommand{\eeg}{\end{exam}}

\newcommand{\cO}{\mathcal{O}}

\newcommand{\ci}{\mathcal{I}}
\def\bz{\mathbb{Z}}

\def\cm{\mathcal{M}}
\def\ci{\mathcal{I}}

\newcommand{\wt}{\widetilde}

\newcommand{\spec}{\textrm{Spec}}

\newcommand{\bff}{\mathbb{F}}

\def\ni{\noindent}

\def\so{\textrm{SO}}
\def\sl{\textrm{SL}}
\def\gl{\textrm{GL}}

\renewcommand{\phi}{\varphi}
\begin{document}
\title[Frobenius splitting of certain rings of invariants]{Frobenius splitting of\\ certain rings of invariants}
\author{V. Lakshmibai}
\address{
Dept. of Mathematics,
Northeastern Univ., 
Boston, MA., U.S.A}
\email{lakshmibai@neu.edu}
\author{K. N. Raghavan}
\address{
The Institute of Mathematical Sciences, 
Chennai 600113, India.}
\email{knr@imsc.res.in}
\author{P. Sankaran}
\address{The Institute of Mathematical Sciences, 
Chennai  600113, India.}
\email{sankaran@imsc.res.in}
\date{}
\thispagestyle{empty}
\footnote{
{\bf keywords:} Frobenius splitting, invariant rings}

\thanks{V.~Lakshmibai was partially supported
by NSF grant DMS-0652386 and Northeastern University RSDF 07-08.}
\thanks{K.~N.~Raghavan and P.~Sankaran were Partially supported
by DAE grant No.~11--R\&D--IMS--5.01--0500.}

\dedicatory{Dedicated to Professor Melvin Hochster on the occasion
of his sixty-fifth birthday.} \maketitle \noindent {\bf Abstract:}
Let $k$ be an algebraically closed field of characteristic $p>0$,
and $V$ an $n$-dimensional $k$-vector space together with a
non-degenerate symmetric bilinear form. Let $G$ denote one of the
groups $G=\textrm{SL}(V)$ or $\textrm{SO}(V)$ where we assume that
$p>2$ if $G=\so(V)$. Let $R$ denote the coordinate ring of
$V_{m,q}:=V^{\oplus m}\oplus V^{*\oplus q}$ (resp. $V_m:=V^{\oplus
m}$) if $G=\textrm{SL}(V)$ (resp. if $G=\textrm{SO}(V)$), $V^*$
being the dual of $V$. The defining representation of $G$ on $V$
induces the diagonal action of $G$ on $V_{m,q}$ (resp. $V_m$). Let
$S=R^{G}$. In this paper, we show that $S$ is Frobenius split.

\noindent
\section{Introduction}

The concept of $F$-purity was introduced by Hochster-Roberts
\cite{hr0}; the $F$-purity for a noetherian ring of prime
characteristic is equivalent to the splitting of the Frobenius
map, when the ring is finitely generated over its subring of
$p$-th powers. It is closely related to the Frobenius splitting
property \'a la Mehta-Ramanathan \cite{mr} for algebraic
varieties; to make it more precise, the $F$-split property for an
irreducible projective variety $X$ over an algebraically closed
field of positive characteristic is equivalent to the $F$-purity
of the ring $\oplus_{n\geq 0} H^0(X;L^n)$ for any ample line
bundle $L$ over $X$ (cf.\cite{hw},\cite{s1},\cite{s2}). We feel
that it is but appropriate to dedicate this paper to Professor
Hochster
 on the occasion of his sixty-fifth birthday and thus make a
modest contribution to this birthday volume.

Let $k$ be an algebraically closed field of characteristic $p>0$
and let $X$ be a $k$-scheme. One has the Frobenius morphism
(which is only an $\bff_p$-morphism) $F\colon X\lr X$
defined as the identity map of the underlying
topological space of $X$, the morphism  of structure sheaves
$F^\#\colon\cO_X\lr \cO_X$
being the $p$-th power map.  The morphism $F$ induces a
morphism of $\cO_X$-modules $\cO_X\lr F_*\cO_X$. The variety $X$ is called
{\it Frobenius split} (or {\it $F$-split} or, merely, {\it split})
if there exists a
splitting $\phi\colon F_*\cO_X\lr \cO_X$ of the morphism $\cO_X\lr F_*\cO_X$.
Equivalently, $X$ is
Frobenius split if there exists a morphism of sheaves of abelian groups
$\phi\colon \cO_X\lr \cO_X$ such that (i) $ \phi(f^pg)=f\phi(g), ~f,g\in \cO_X$
and (ii) $\phi(1)=1$.
Basic examples of varieties that are Frobenius
split are smooth affine varieties, toric varieties (cf. \cite{bk}),
generalized flag varieties, and Schubert varieties \cite{mr}. Smooth
projective curves of genus greater than $1$ are examples
of varieties that are {\it not} Frobenius split.

Frobenius splitting is an interesting property to study. If $X$ is
Frobenius split, then it is weakly normal (cf. \cite{bk}, Prop
1.2.5) and reduced (cf. \cite{bk}, Prop. 1.2.1). Indeed,
projective varieties which are Frobenius split are very special.
We refer the reader to \cite{bk} for further details.

If $X=\textrm{Spec}(R)$, then $X$ is Frobenius split if and only
if the Frobenius homomorphism $R\lr R$ defined as $a\mapsto a^p$
admits a splitting $\varphi\colon R\lr R$ such that $\varphi(a^pb)=a\varphi(b)$,
and $\phi(1)=1$.

If a linearly reductive group $G$ acts on a $k$-algebra $R$ which
is Frobenius split, then the invariant ring $R^G$ is Frobenius
split  (see \cite[Exercise 1.1..E(5)]{bk}). To quote Karen Smith
\cite[p. 571]{s1}, ``The story of $F$-splitting and global
$F$-regularity for quotients by reductive groups in
characteristics $p$ that are not linearly reductive is much more
subtle and complicated". We shall show that although the groups
$\so(n), n\geq 3,$ and $\sl(n), n\geq 2,$ are not linearly
reductive, it turns out that certain rings of invariants for these
groups {\it are} Frobenius split.

We state below the main results of this paper.

Let $k$ be an algebraically closed field of characteristic $p>2$
and $V$ an $n$-dimensional vector space over $k$ with a
symmetric non-degenerate bilinear form. Denote by $A_m$ the
the coordinate ring of $V_m:=V^{\oplus m}, m\geq 1,$ and consider the
action of $\so(V)$ on $A_m$ induced by the diagonal action of
$\so(V)$ on $V^{\oplus m}$. Then

\bth\label{soninvariants}
The invariant ring $A_m^{\so(V)}$ is Frobenius split for all $m\geq 1$.
\eeth

The group $\sl(V)$ acts on $V$, as well as on the dual vector
space $V^*=\textrm{Hom}_k(V,k)$. Now consider the diagonal action
of $\sl(V)$ on the vector space $V_{m,q}:=V^{\oplus m}\oplus
V^{*\oplus q}$, $m,q\geq 1$. This leads to an action of $\sl(V)$
on the coordinate ring $A_{m,q}$ of $V_{m,q}$.

\bth\label{slninvariants} The invariant  ring $A_{m,q}^{\sl(V)}$
is Frobenius split for any $m,q\geq n$. \eeth

We shall now sketch the proofs of the main results (assuming $m, q
>n$). Let $S$ be the invariant ring in question. Let $R$ be the ring of invariants
under the larger group $\wt{G}=\gl(V)$  (resp.
$\wt{G}=\textrm{O}(V)$), we have (cf. \cite{d-p,ls1}) that $R$ is
the coordinate ring  of a certain determinantal variety  in
$M_{m,q}$, the space of $m\times q$ matrices (resp. $Sym\,M_m$,
the space of symmetric $m\times m$ matrices) with entries in $k$.
Now a determinantal variety  in $M_{m,q}$ (resp. $Sym\,M_m$) can
be canonically identified (cf. \cite{ls1}) with an open subset in
a certain Schubert variety in $G_{q,m+q}$, the Grassmannian
variety of $q$-dimensional subspaces of $k^{m+q}$ (resp. the
symplectic Grassmannian variety, the variety of all maximal
isotropic subspaces of a $2m$-dimensional vector space over $k$
endowed with a non-degenerate skew-symmetric bilinear form). Hence
we obtain that $R$ is Frobenius split (since Schubert varieties
are Frobenius split). Let  $X=\spec(S),Y=\spec(R) $, and
$\pi\colon X\lr Y$,  the morphism induced by the inclusion
$R\subset S$ . When $G=\so(V)$, we show that $\pi$ is a double
cover which is \'etale over a `large' open subvariety -- that is a
subvariety whose complement has codimension at least $2$. When
$G=\sl(n)$, we show that restricted to a large open subvariety,
$\pi$ is a $\mathbb{G}_m$ bundle. The main theorems are then
deduced using normality of $S$. 

Theorem \ref{soninvariants} can also be deduced from Hashimoto's
work~\cite{hash}, wherein it is shown that if a reductive group
$G$ acts on a polynomial ring $A$ over $k$ (of positive
characteristic) with good filtration, then the ring $A^G$ of
invariants is strongly $F$-regular. Our Theorem
\ref{slninvariants} does not seem to follow from the results of
\cite{hash}.  Granting the results of \cite{ls} and
\cite{lrss}---we don't need all the results of these papers, only
some of the relatively easier ones---the arguments used in our
proofs are straightforward and quite elementary; the techniques
used in~\cite{hash} are representation theoretic.

Theorem \ref{soninvariants} will be proved in \S 2 and
Theorem \ref{slninvariants} in \S 3.

\section{Splitting $\so(n)$-invariants}

Let $k$ be an algebraically closed field of characteristic $p>0$.
Suppose that $S$ is an affine $k$-algebra which is Frobenius split and that
a finite group $\Gamma$ acts on $S$ as $k$-algebra automorphisms.
Then the invariant ring $R=S^{\Gamma}$ is Frobenius split provided
the order of $\Gamma$ is not divisible by $p$ (cf. \cite[Ex. 1.1.E(5)]{bk}).
We first obtain a partial converse to this statement in the case when
$\Gamma$ is of order $2$.

Assume that $\textrm{char}(k)>2$. Let $S$ be an affine $k$-domain and let
$\Gamma=\{1,\gamma\}\cong \bz/2\bz$ act effectively on $S$. Denote by  $R$ the
invariant subalgebra $S^\Gamma$.   Then~$R$ is an affine $k$-algebra and~$S$ is quadratic and integral over $R$. Indeed,
any $b\in S$ can be expressed as $b=b_0+b_1$ where $b_0=(1/2) (b+\gamma(b))
\in R$ and $b_1=(1/2)(b-\gamma(b))$ satisfies $\gamma(b_1)=-b_1$. Thus,
we can choose generators $u_1,\cdots, u_s$ for the $R$-algebra
$S$ to be in the $-1$ eigenspace of $\gamma$. Clearly $u_i^2=-u_i\gamma(u_i)
=:p_i\in R$ for all $i\leq s$.  Furthermore,

\ni $\gamma(u_iu_j)=u_iu_j=:p_{i,j}\in R$ for all $i,j\leq s$
(with $p_{i,i}=p_i$). Observe that $p_{i,j}^2=p_ip_j$.

We shall assume that $S$ is reduced so that $p_i\neq 0$, for all
$i$. Now let $R_i=R[1/p_i], 1\leq i\leq n$. Let
$S_i=S[1/u_i]$. We claim that $S_i=R_i[u_i]$. To see this, first
observe that $R_i[u_i]\subset S[1/u_i]$, since $1/p_i=(1/u_i)^2\in
S[1/u_i]$. To show that $S[1/u_i]\subset R_i[u_i]$, it suffices to
show that $u_j\in R_i[u_i]$ for all $j$ and $(1/u_i)\in R_i[u_i]$.
Indeed, $1/u_i=u_i/u_i^2=u_i/p_i\in R_i[u_i]$ and so
$u_j=p_{i,j}/u_i\in R_i[u_i]$.

Write $X=\textrm{Spec}(S), Y=\textrm{Spec}(R)$ and let $\pi:X\lr
Y$ be the morphism (induced by the inclusion $R\subset S$). As
above, let $S_i=S[1/u_i]$, and let $U_i=\spec(S_i)\subset X$ and
let $U:=\bigcup_{1\leq i\leq s} U_i$; it is the full inverse image
under $\pi$ of $\bigcup_{1\leq i\leq s} \spec(R_i)$. It is readily
verified that $\pi|U \colon U\lr \pi(U)$ is \'etale. Indeed, $S_i$
is a free $R_i$ module with basis $\{1,u_i\}$ and
$det(u_i)=-p_i\neq 0$ and so $\pi|U_i$ is \'etale. Hence $\pi|U$
is \'etale.

On the other hand, if $y\in Y$ is a closed point
such that $p_i(y)=0$ for all $i\leq s$,
then the fibre $f^{-1}(y)=\spec(S_y\otimes_{R_y}k)$ is the
scheme whose coordinate ring is
$S_y\otimes_{R_y}k=k[u_1,\cdots, u_s]/\langle u_i^2, ~1\leq i\leq s\rangle$. Here
$R_y$ is the local ring at $y$. Thus  $f^{-1}(y)$ is
non-reduced. It follows
that the ramification locus of $\pi$ {\it equals} $Y\setminus \pi(U)$.
(See  \cite[\S III.10, Theorem 3]{mum}.)

\noindent
\bpropo \label{affinecover}
Let $k$ be an algebraically closed field of characteristic $p>2$.
Let $S$ be an affine normal domain over $k$ acted on by $\Gamma\cong \bz/2\bz$
and let $R:=S^{\Gamma}$ be Frobenius split. Suppose
that the ramification locus of the double cover
$\pi\colon \textrm{Spec}(S)\lr \textrm{Spec}(R)$
has codimension at least $2$ in $\textrm{Spec}(R)$.
Then, any splitting $\phi:R\lr R$
extends uniquely to a splitting $\psi\colon S\lr S$.

\epropo
\noindent
{\it Proof.} We use the notations introduced above.

Since $X$ is normal and the codimension of the the ramification
locus of $\pi$ is at least $2$, it suffices to show that
$U=\bigcup_{1\leq i\leq s} U_i$ is Frobenius split (cf. \cite{bk},
Lemma 1.1.7,(iii)).

Let $\phi:R\lr R$ be a splitting of $Y=\spec(R)$. First, we shall
 extend $\phi$ to a splitting $\psi_i\colon S_i\lr S_i$ of
$U_i=\spec(S_i)(=\spec(R_i[u_i]))$ for each $i$ and verify that
these splittings agree on the overlaps $U_i\cap U_j$ for $1\leq
i,j\leq s$. Thus we will obtain a splitting of $U=\bigcup_{1\leq
i\leq s} U_i$. By normality of $X$ and the hypothesis on the
codimension of the ramification locus, we will conclude that this
splitting extends to a splitting of $X$. Next, we shall
 establish the uniqueness of the extension.

Recall that  $\{1, u_i\}$ is an $R_i$-basis for $S_i$.  Since
$u_i=u_i^{-p}p_i^{(1+p)/2}$ on $U_i$,  if $\psi_i\colon S_i\lr
S_i$ is {\it any} splitting of $U_i$ which extends the splitting
$\phi_i$ of $R_i$ defined by $\phi$, we must have
$\psi_i(au_i)=\psi_i((1/u_i)^pap_i^{(p+1)/2})
=(1/u_i)\phi_i(ap_i^{(p+1)/2})$. By additivity, we must have
$$\psi_i(au_i+b)=(1/u_i)\phi_i(ap_i^{(p+1)/2})+\phi_i(b)
=(u_i/p_i)\phi_i(ap_i^{(p+1)/2})+\phi_i(b)$$ where $a,b\in R_i$.
Thus the extension $\psi_i$, if it exists, is unique.

We now {\it define} $\psi_i$ by the above equation and claim that
$\psi_i$ is indeed a splitting of $S_i$. First, observe that
$\psi_i(1)=1$, by the very definition of $\psi_i$.

Now, for any $x,y,a\in R_i$, we have

\ni $\psi_i((xu_i+y)^pau_i)=\psi_i(x^pp_i^{(p+1)/2}a+y^pau_i)$

\ni $=x\phi(p_i^{(p+1)/2}a)+y\phi(au_i)$

\ni $=x(p_i/u_i)\psi_i(au_i)+y\phi(au_i)$

\ni $=xu_i \psi_i(au_i)+y\psi(au_i)$

\ni $=(xu_i+y)\phi(au_i)$

An entirely similar (and easier) computation shows that

\ni $\psi_i((xu_i+y)^pb)=(xu_i+y)\psi_i(b)$, completing the
verification that $\psi_i$ is a splitting.

We verify, by another straightforward computation, that these $\psi_i$
agree on the overlaps $U_i\cap U_j$.
Indeed, writing $u_j=u_ip_{i,j}/p_i,$ we have \\
$\psi_i(u_j)=\psi_i(u_ip_{i,j}/p_i)
=(u_i/p_i)\phi((p_{i,j}/p_i)p_i^{(p+1)/2})$

\ni $=(u_i/p_i)\phi(p_{i,j}p_i^{(p-1)/2})$

\ni Since $p_i=p_{i,j}^2/p_j $ on $U_i\cap U_j$, we have

\ni $\phi(p_{i,j}p_i^{(p-1)/2})
=\phi(p_{i,j}^p(p_i/p_{i,j}^2)^{(p-1)/2})$

\ni $=p_{i,j}\phi(p_j^{(1-p)/2})=(p_{i,j}/p_j)\phi(p_j^{(p+1)/2})$

\ni Substituting in the above expression for $\psi_i(u_j)$ we get
$$\psi_i(u_j)=(u_ip_{i,j}/(p_ip_j))\phi(p_j^{(p+1)/2})
=(u_j/p_j)\phi(p_j^{(p+1)/2})=\psi_j(u_j)$$

This implies that the extensions $\{\psi_i\}$ patch to yield
a well-defined splitting on $U$ as claimed. As observed above, the
normality of $X$ and the hypothesis on the codimension of the
ramification locus implies that
this splitting extends to a {\it unique} splitting $\psi\colon S\lr S$.

Finally, if $\psi'$ is another splitting of $X$ which also extends $\phi$,
then $\psi'$ and $\psi$ agree on $U_i$ (for any $i$) as already observed.
As $X$ is irreducible, $U_i$ is dense in $X$ and
we conclude that $\psi'=\psi$. \hfill $\Box$

As a corollary, we obtain the following

\bth\label{schemecover}
Let $\pi:X\lr Y$ be a double cover of a Noetherian scheme whose
ramification locus has codimension at least $2$. Suppose
that $X$ is reduced, irreducible and normal and that $Y$ is Frobenius split,
then $X$ is Frobenius split.
\eeth

\noindent
{\it Proof.}  Cover $X$ by finitely many affine patches $X_\alpha$.
Let $Y_\alpha:=\pi X_\alpha$.
Then each $\pi|_{X_\alpha}$ satisfies the hypotheses of the above
proposition. Let~$\phi$ be a splitting of $Y$ and let
$\psi_\alpha$ be the unique splitting of $X_\alpha$ that extends the
splitting $\phi|_{Y_\alpha}$.
The $\psi_\alpha$'s agree on overlaps and hence define a unique
splitting of $X$ which `extends' $\phi$. \hfill $\Box$

We now turn to proof of Theorem \ref{soninvariants}.

\noindent {\it Proof of Theorem \ref{soninvariants}.} Denote by
$S$ the ring of $\so(V)$-invariants of $A_m,$ where $A_m$ is the
coordinate ring of $V_m$.   Let~$R$ be the ring of $\textrm{O}(V)$-invariants.

We shall  assume that $m>n$.
  By \cite{d-p,ls1} we have  that
$Y:=\spec(R)$ is the determinantal variety $D_n(Sym\,M_m)$
consisting of all matrices in $Sym\,M_m$ (the space of symmetric
$m\times m$ matrices with entries in $k$) of rank at most $n$;
further, we have (cf. \cite{ls1}) an identification of
$D_n(Sym\,M_m)$  with an open subset of a certain Schubert variety
in the Lagrangian Grassmannian variety (of all maximal isotropic
subspaces of a $2m$-dimensional vector space over $k$ endowed with
a non-degenerate skew-symmetric bilinear form). Hence we obtain
that $Y$ is Frobenius split (since Schubert varieties are
Frobenius split (cf. \cite{mr})).

Observe that $\Gamma:=\textrm{O}(n)/\so(n)\cong \bz/2\bz$ acts on
$S$ (the subring of $\textrm{SO}(V)$-invariants of $A_m$) and that
$S^\Gamma$ equals $R$. As above, let

\ni $X:=\spec(S)$, and $\pi:X\lr Y$ be the morphism induced by the
inclusion $R\subset S$. We need only verify the hypotheses of
Theorem \ref{schemecover}. It is well-known that $S$ is an affine
normal domain. It remains to verify that the codimension of the
branch locus is at least two. This was proved in \cite{lrss}. In
fact, the ramification locus of $Y$ {\it equals} the singular
locus of $Y$, but this more refined assertion is not relevant
here. Since $Y$ is normal it follows that the codimension of the
ramification locus is at least $2$.  Therefore, by Theorem
\ref{schemecover}, $X$ is Frobenius split.

The case $m=n$ is isolated separately as Lemma~\ref{mequalsn} below.
When $m<n$, it is easy to see that $S=R$. Again, $R$ is
a polynomial algebra over $k$ and hence $S$ is Frobenius split. \hfill $\Box$

Assume that $m=n$.  
In this case $R=k[y_{i,j}:1\leq i\leq j\leq n]$ is a polynomial
ring, being the ring of polynomial functions on the space of
$n\times n$ symmetric matrices.    As an $R$-algebra,
$S=R[u]/\langle u^2-f\rangle$, where $f$ denotes the determinant
function of the symmetric $n\times n$ matrix whose entry in
position $(i,j)$ for $1\leq i\leq j\leq n$ is $y_{i,j}$. \noindent
\blem\label{mequalsn} Let $m=n$.     The ring~$S$ of ${\rm
SO}(V)$-invariants is Frobenius split in this case also. \elem
\noindent {\it Proof.} There is a natural identification of
$\spec(R)$ with an affine patch of the symplectic Grassmannian and
the vanishing locus of~$f$ under this identification becomes an
open part of a Schubert variety~\cite{ls1,lrss}.   Thus
by~\cite{mr} (see also~\cite{bk}), there exists a splitting of
$\spec(R)$ that compatibly splits~$\spec(R/(f))$. Let~$\phi$ be
such a splitting.   Continue to denote by~$\phi$ the restriction
of~$\phi$ to the open part $\spec(R[1/f])$.   Arguing as in the
proof of Proposition~\ref{affinecover} above,  we may `lift' the
restriction~$\phi$ to a splitting (also denoted~$\phi$) of
$\spec(S[1/f])$.    We claim that $\phi$ maps~$S$ to~$S$ and hence
extends to a splitting of~$\spec(S)$.   Indeed,  a general element
of~$S$ is of the form~$au+b$ with $a$,~$b$ in~$R$,   so that
$\phi(au+b)=\phi({\frac{au^{p+1}}{u^p}}\ +\
b)={\frac{\phi(af^{(p+1)/2})}{u}}\ +\ \phi(b)$. Since~$\phi$
compatibly splits the vanishing locus of~$f$, it follows that
$\phi(af^{(p+1)/2})$ belongs to the ideal~$(f)$. Writing
$\phi(af^{(p+1)/2})=cf$,  we have
$\phi(au+b)={\frac{cf}{u}}+\phi(b)=cu+\phi(b)\in S$. \hfill $\Box$

We conclude this section with the following remarks.

\noindent
\brem (i) The condition on codimension of $U$ in
Proposition \ref{affinecover}  will be satisfied
if $S$ is generated over $R$ by two or more elements~$u_i$ such that
there exist $u_i, u_j$ such that the supports $D_i$ and $D_j$ of the
reduced scheme defined by $u_i=0$ and $u_j=0$ do
not have any component in common.

\noindent (ii) Theorem \ref{schemecover} is not valid when the
hypothesis on the codimension of the ramification locus is not
satisfied. For example, if $\Gamma\cong \bz/2\bz$ is generated by
the involution of a hyperelliptic curve $X$ of genus $g\geq 2$,
then the quotient is a smooth projective curve which is Frobenius
split. However, $X$ is not split since $\omega_X$ is ample but
$H^1(X;\omega) \cong k$, whereas higher cohomologies for ample
line bundles over Frobenius split projective varieties vanish.

\noindent
(iii) We do not know if Theorem \ref{schemecover} remains valid
if $\Gamma$ is {\it any} finite group whose order is prime to the
characteristic $p$ of $k$, even in the
case when $\Gamma$ is cyclic.

\noindent
(iv) One has an isomorphism of $\sl(2)$ with $\so(3)$ such that
the $\so(3)$ action on $V=k^3$ corresponds to the conjugation
action of $\sl(2)$ on trace zero $2\times 2$ matrices. In this
case the Frobenius splitting property of $A_m$ was proved by
Mehta-Ramadas \cite[Theroem 6]{mr2}.  It should be noted that when
$\dim V=3$, the completion of the stalk at the origin in $A_m$ is
isomorphic to the completion of the stalk at the point
corresponding to the class of the trivial rank $2$ vector bundle
in the moduli space of equivalence classes of semi-stable, rank
$2$ vector bundles with trivial determinant on a smooth projective curve of
genus $m> 2$ (see \cite{mr2}). \erem

\section{Splitting $\sl(n)$ invariants}

In this section we shall establish Theorem \ref{slninvariants}.
Let $V$ be an $n$ dimensional vector space over an algebraically
closed field $k$ of characteristic $p\geq 2$ and let $V^*$ denote
its dual.  Let $V_{m,q}:=V^{\oplus m}\oplus V^{*\oplus q}$, and
let $A$ denote the ring of regular functions on $V_{m,q}$. By
fixing dual bases for $V$ and $ V^*$, we shall view elements of
$V$ and $V^*$ as row and column vectors respectively, so that
$V^{\oplus m}$ (resp. $V^{*\oplus q}$) is identified with the
space $M_{m,n}$ of $m\times n$ matrices (resp. the space $M_{n,q}$
of $n\times q$ matrices) over $k$; further, $GL(V)$ gets
identified with $\gl_n(k)$ (the group of invertible $n\times n$
matrices over $k$). In the sequel, we shall
 denote $\gl_n(k)$ by just $\gl(n)$. Then the action of
$\textrm{GL}(V)$ on $V^{\oplus m} $ gets identified with the
multiplication on the right of $M_{m,n}$ by $\textrm{GL}(n)$.
Similarly, the action of $g\in \gl(V)$ on $V^{*\oplus q}$ gets
identified with the multiplication on the left of $M_{n,q}$ by
$g^{-1}$. The diagonal action of $\gl(V)$ on $V^{\oplus m}\oplus
V^{*\oplus}$ is therefore defined as $(u,\xi).g=(ug,g^{-1}\xi)$
where $g\in \gl(n)$ and $(u,\xi)\in \cm_{m,q}:=M_{m,n}\times
M_{n,q}$. We identify $A$ with the coordinate ring of $\cm_{m,q}$.

We denote by $R$ and $S$ the
rings of invariants $A^{\textrm{GL}(n)}$ and $A^{\sl(n)}$ respectively.
Let $Y=\spec(R)$ and $X=\spec(S)$.  Note that $Y$ and $X$ are the
GIT quotients $\cm_{m,q}/\!\!/\textrm{GL}(n)$ and $\cm_{m,q}
/\!\!/\textrm{SL}(n)$ respectively.

Let $m,q\geq n$. By \cite{d-p,ls1} we have  that $Y$ is the
determinantal variety $D_n(M_{m,q})$ consisting of all matrices in
$M_{m,q}$ (the space of $m\times q$ matrices with entries in $k$)
of rank at most $n$; further, we have (cf. \cite{ls1}) an
identification of $D_n(M_{m,q})$  with an open subset of a certain
Schubert variety in the Grassmannian variety (of $q$-dimensional
subspaces of $k^{m+q}$). Hence we obtain that $Y$ is Frobenius
split (since Schubert varieties are Frobenius split). The
multiplication map $\mu\colon \cm_{m,q}\,\lr \,M_{m,q}$ factors
through $Y$; further, under $\pi\colon X\lr Y$ (induced by the
inclusion $R\subset S$), we have, $\pi([u,\xi])=u\xi\in M_{m,q}$
where $[u,\xi]\in X$ is the image of $(u,\xi)\in \cm_{m,q}$ under
the GIT quotient $\cm_{m,q}\lr X$.

 Let $I(n,m)$ denote the set of all $n$-element
subsets $I$ of $\{1,2,\cdots,m\}$. Any such $I$ determines a
regular function $u_I\colon \cm_{m,q}\,\lr\, k$ which maps
$(u,\xi)$ to the determinant of the $n\times n$ submatrix $u(I)$
of $u\in M_{m,n}$ with column entries given by $I$.
Clearly $u_I$ is invariant under the action
of $\sl(n)$ on $\cm_{m,q}$ and hence yields a regular function $u_I$ on $X$.

We define $\xi(J)$ and $\xi_J$ for $J\in I(n,q)$ analogously;
$\xi_J$ is also an  $\sl(n)$-invariant.

We have, $u_I\xi_J=:p_{I,J}\in R$ for all $I\in I(n,m), J\in
I(n,q)$; indeed $p_{I,J}([u,\xi])$ is just the determinant of the
$n\times n$ submatrix of $u\xi\in M_{m,q}$ with row and column
indices given by $I$ and $J$ respectively. It is shown in
\cite{ls}, among other things, that $S$ is generated as an
$R$-algebra by $u_I,\xi_J, ~I\in I(n,m),J\in I(n,q)$, the ideal of
relations being generated by $u_I\xi_J-p_{I,J}, ~I\in I(n,m),J\in
I(n,q)$ together with certain quadratic relations among the
$u_I$'s and certain quadratic relations among the $\xi_J$'s.
Further, in \cite{ls}, a standard monomial basis is constructed
for $S$; as a particular consequence, we have that each $u_I$
(resp. $\xi_J$) is algebraically independent over $R$ for $I\in
I(n,m)$ (resp. $J\in I(n,q)$).

\ni For $K\in I(n,m), L\in I(n,q)$, let
$$R_{K,L}=R[1/p_{K,L}],Y_{K,L}=\spec(R_{K,L})$$ For a given $I\in I(n,m), J\in I(n,q)$, let
$$Y_I={\underset{J'\in I(n,q)}{\cup}}Y_{I,J'},\ Y_J={\underset{I'\in I(n,m)}{\cup}}Y_{I',J}$$
 Note that for $I\in I(n,m)$, any $Y_{I,J'}$ is contained in $Y_I$; similarly, for $J\in
I(n,q)$, any $Y_{I',J}$ is contained in $Y_J$.

Set $X_I=\pi^{-1}(Y_I)\subset X$ and $X_J=\pi^{-1}(Y_J)\subset X$.
Note that $u_I$ (resp. $\xi_J$) is non-zero on $X_I$ (resp.
$X_J$). Denote by $f_I\colon X_I\lr Y_I\times k^*$ the morphism
$f_I=(\pi|X_I,u_I|X_I)$, and by $f_J\colon X_J\lr Y_J\times k^*$
the morphism $f_J=(\pi|X_J,\xi_J|X_J)$.

\noindent \blem \label{localtriviality} The morphisms $f_I\colon
X_I\lr Y_I\times k^*$ and $f_J\colon X_J\lr Y_J\times k^*$ are
isomorphisms for any $I\in I(n,m), J\in I(n,q)$. \elem \noindent
{\it Proof.}  We shall prove that $f_I$ is an isomorphism, the
proof in the case of $f_J$ being the same.

 Let $X_{I,J}=\pi^{-1}(Y_{I,J})$; then $X_{I,J}$ equals
$\spec(S_{I,J})$ (where $S_{I,J}=S[1/p_{I,J}]$) and $X_{I,J}$ is
contained in $X_I$. The morphism $f_{I,J}\colon X_{I,J}\lr
Y_{I,J}$ defined by the restriction of $f_I$ is induced by the
$R_{I,J}$-algebra map $f_{I,J}^*\colon R_{I,J}[t,t^{-1}]\lr
S_{I,J}$ which maps $t$ to $u_{I}$. Note that $p_{I,J}=u_I\xi_J$
implies that $u_I$ is invertible in $S_{I,J}(=S[1/p_{I,J}])$.

We must show that
\begin{enumerate}
\item \label{iso} $f_{I,J}^*$ is an isomorphism of $k$-algebras
\item \label{patch} $f_{I,J}$ and $f_{I,J'}$ agree on the overlap
$X_{I,J}\cap X_{I,J'}$ for any two $J,J'\in I(n,q)$.
\end{enumerate}

\noindent (\ref{iso}).  Note that
$u_{I'}=u_Iu_{I'}\xi_J/p_{I,J}=u_Ip_{I',J}/p_{I,J}=f_{I,J}^*(p_{I',J}/p_{I,J}t)$.
Hence $u_{I'}$ is in the image of $f_{I,J}^*$ for any $I'\in
I(n,m)$. Similarly $\xi_{J'}$ is also in the image of $f_{I,J}^*$
for any $J'\in I(n,q)$. Therefore $f_{I,J}^*$ is surjective. Now
suppose that $P(t)\in R_{I,J}[t,t^{-1}]$ is in the kernel of
$f^*_{I,J}.$ We may assume that $P(t)$ is a polynomial in $t$ and
that the coefficients of $P(t)$ are actually in $R$. Then
$0=f_{I,J}^*(P(t))=P(u_I).$  Since $X_{I,J}$ is open in $X$, which
is irreducible, the we see that the equation $P(u_I)=0$ must hold
in $S$. This contradicts the fact that $u_I$ is algebraically
independent over $R$ (cf. \cite{ls},Theorem 6.06,(3)).  Hence
$f_{I,J}^*$ is an isomorphism.

\noindent
(\ref{patch}). It is evident that
$f_{I,J}^*(t)=u_I\in S_{I,J}$ and $f_{I,J'}^*(t)=u_I\in
S_{I,J'}$ both restrict to the same regular function, namely
$u_I|X_{I,J}\cap X_{I,J'}$, on the overlap $X_{I,J}\cap X_{I,J'}=
\spec(S[1/p_{I,J},1/p_{I,J'}])$.  It follows that $f_{I,J}$ and $f_{I,J'}$
agree on $X_{I,J}\cap X_{I,J'}$.  This completes the proof that
$f_I$ is an isomorphism. \hfill $\Box$

Observe that, if $J,J'\in I(n,q)$,
then $\xi_J/\xi_{J'}\in S[1/\xi_{J'}]$ defines a regular function on $Y_{J'}$.
This is because, on $Y_{I,J'}$, $\xi_J/\xi_{J'}=(u_{I}\xi_J)/(u_I\xi_{J'})
=p_{I,J}/p_{I,J'}$. It is immediately seen that, on $Y_{I,J'}\cap Y_{I',J'},$
the two regular functions $p_{I,J}/p_{I,J'}$ and $p_{I',J}/p_{I',J'}$
agree.  Therefore we conclude that $\xi_J/\xi_{J'}$ is a well-defined
regular function on $Y_{J'}$.  Clearly it is invertible on $Y_J\cap Y_{J'}$.
Similar statements concerning $u_I/u_{I'}$ hold for any $I,I'\in I(n,m)$.

\noindent
{\bf  Notation:}  Let $m,q\geq n$.

\ni Denote by $\ci$ the disjoint union $I(n,m)\coprod I(n,q)$.  We
set
$$\lambda_{\beta,\alpha}=\left\{
\begin{array}{ccl}u_\alpha/u_\beta &
~if~& \alpha,\beta\in I(n,m), \\
\xi_\beta/\xi_\alpha &~if ~&\alpha,\beta\in I(n,q),\\
p_{\alpha,\beta} &~if ~&\beta\in I(n,q),\alpha\in I(n,m),\\
1/p_{\beta,\alpha}&~if~&\beta\in I(n,m),\alpha\in I(n,q).\\
\end{array}\right.$$

\ni Consider the covering $\{Y_\alpha\}_{\alpha\in \ci}$
of the open subvariety $Y_0:=\bigcup_{\alpha\in \ci}Y_{\alpha}\subset Y$.
The cocycle condition
$\lambda_{\alpha,\beta}\lambda_{\beta,\gamma}= \lambda_{\alpha,
\gamma}$ is readily verified for any $\alpha, \beta,\gamma \in
\ci$. Thus we obtain a $\mathbb{G}_m$-bundle over $Y_0$; call it
$\mathcal{E}$.

Let $X_0:=\bigcup_{\alpha\in \ci} X_{\alpha}$.

\noindent \blem \label{linebundle}  Assume that $m,q\geq n$. With
the above notations, the total space of the $\mathbb{G}_m$-bundle
$\mathcal{E}$ over $Y_0$ is isomorphic to the open subvariety
$X_0:=\bigcup_{\alpha\in\ci}X_\alpha\subset X$. \elem

\noindent
{\it Proof.} The total space of the ${\Bbb G}_m$-bundle corresponding to
$\mathcal{D}$ is \newline $\coprod_{\alpha\in \ci }
Y_\alpha\times k^*/\!\!\sim$ where
$(\pi([u,\xi]), t)\in Y_\alpha\times k^*$ is identified with \newline
$(\pi([u;\xi]),
\lambda_{\beta,\alpha}(\pi([u,\xi]).t))\in Y_\beta\times k^*$ whenever
$\pi([u,\xi])\in Y_\alpha\cap Y_\beta$.
One has the following commuting diagram for any
$\alpha,\beta\in\ci$:
$$\begin{array}{cccccccc}
Y_\alpha\times k^*&\supset& (Y_\alpha\cap Y_\beta)\times k^*&
\stackrel{\lambda_{\beta,\alpha}}{\lr}&
(Y_{\beta}\cap Y_{\alpha})
\times k^*& \subset &Y_\beta\times k^*\\
f_\alpha \uparrow && f'_\alpha\uparrow && \uparrow f'_\beta &&
$$ \uparrow f_\beta\\
~X_\alpha &\supset&  X_\alpha\cap X_\beta & == &X_\beta\cap X_\alpha &\subset
& X_\beta~\\
\end{array}$$
where $f'_\alpha$ is the restriction of $f_\alpha$. Since, by
Lemma \ref{localtriviality}, the $f_\alpha$'s are isomorphism of
varieties, it follows that that the total space of the
$\mathbb{G}_m$ bundle over $Y_0$ is isomorphic to the union
$X_0:=\bigcup_\alpha X_\alpha\subset X$. \hfill $\Box$

We shall now compute the codimension of $Z:=X-X_0$. We give the
reduced scheme structure on $Z$.
It is evident that $Z$ is defined by the
equations $p_{I,J}=0, \forall I\in I(n,m),J\in I(n,q).$
We claim
$Z=Z_u\cup Z_\xi$ where $Z_u$ is the closed
subvariety with reduced scheme structure
defined by the equations $u_I=0, \forall I\in I(n,m)$ and $Z_\xi$,
by the equations $\xi_J=0,\forall J\in I(n,q)$.
Clearly $Z_u\cup Z_\xi\subset Z$.
On the other hand, if $[u,\xi]$
is not in $Z_u\cup Z_\xi $, then $u_I([u,\xi])\neq 0$ for some $I$ and
$\xi_J([u,\xi])\neq 0$ for some $J$. This implies that $p_{I,J}([u,\xi]
\neq 0$. Hence $[u,\xi]\in X_0$.
Thus $Z_u\cup Z_\xi=Z$.

\noindent \blem \label{codimension} Let $m>n$ (resp. $q>n$).  Then
the codimension of $Z_u$ (resp. $Z_\xi$) in $X$ is at least $2$.
\elem \noindent{\it Proof.} Consider the closed subvariety
$M_u:=D_{n-1}(M_{m,n})\times M_{n,q}\subset \cm_{m,q}$ (with
reduced scheme structure). We have,

\ni $dim\,M_u=(n-1)(m+1)+nq$ (note the dimension of the
determinantal variety $D_{t}(M_{r,s})$ (consisting of $r\times s$
matrices of rank at most $t$) equals $t(r+s-t)$ (cf. \cite{ls1})).
Clearly $M_u$ is stable under the $\sl(n)$-action and
$M_u/\!\!/\sl(n)=Z_u$. We shall find an open subset $Z_{u,0}$ of
$Z_u$ such that $\sl(n)$ acts {\it freely} on the inverse image of
$Z_{u,0}$ under the quotient morphism $\eta\colon M_u\lr Z_u$ and
$\eta^{-1}(Z_{u,0})/\!\!/\sl(n) =Z_{u,0}$. It would then follow
that
$\dim(Z_u)=\dim(\eta^{-1}(Z_u))-\dim(\sl(n))=(n-1)(m+1)+nq-(n^2-1)
=(m+n)q-(n^2-1)-(m-n+1)=\dim(X)-(m-n+1)\leq\dim(X)-2$ (note that
$dim\,X=(m+n)q-(n^2-1)$ (cf. \cite{ls})).

Define $$W_u=D_n(M_{m,n})\times M^0_{n,q}$$ where $M^{0}_{n,q}:=
\{\xi\in M_{n,q}\mid\xi_J(\xi)\neq 0,\ \textrm{for some}~J \in
I(n,q)\}.$  Then $W_u$ is the inverse image of
$$Z_{u,0}:=\{[u,\xi]\mid \xi_J(\xi)\neq 0 \}$$ under the quotient
morphism $\eta\colon M_{u}\lr Z_u$. The assertion that the
$\sl(n)$-action is free on $W_u$ follows from the fact that the
$\sl(n)$-action on $M^0_{n,q}$ is free.

An entirely similar argument shows that $Z_\xi$ has codimension at
least $2$, and consequently codimension of $Z$ in $X$ is at least
$2$. \hfill $\Box$

We are now ready to prove Theorem \ref{slninvariants}.

\noindent {\it Proof of Theorem \ref{slninvariants}.} Let $m,q>n$.
As already observed, we have that $Y=D_{n}(M_{m,q})$ can be
identified with an open subset of a certain Schubert variety in
the Grassmann variety $\sl(m+q)/P_q$ of $q$-dimensional vector
subspaces in $k^{m+q}$. Since Schubert varieties in the Grassmann
variety are Frobenius split, it follows that $Y$ is Frobenius
split. Since $Y_0$ is open in $Y$, it follows that it is also
Frobenius split. The variety $X_0$, being the total space of a
$\mathbb{G}_m$-bundle over $Y_0$, is Frobenius split by
\cite{bk},Lemma 1.1.11. Now $X$ being normal and codimension of
$X_0$ in $X$ being at least $2$, it follows that $X$ is Frobenius
split (cf. \cite{bk},Lemma 1.1.7,(iii)).

If $m,q<n$, then $X=Y=M_{m,q}$ and hence Frobenius split.
The case $m=n$ is isolated separately as Lemma~\ref{slmequalsqequalsn} below.
\hfill $\Box$

Assume that $q=n=m$.  In this case $Y=M_{n,n}$.  Denote the
$(i,j)$-th coordinate function on $Y$ by $y_{i,j}$.  The set
$I(n,m)$ and $I(n,q)$ are singletons and so $S=R[u,\xi]/(u\xi-f)$
where $f$ is the determinant function on $Y=M_{m,q}$.
\blem\label{slmequalsqequalsn} Let $q=n=m$.    The ring $S$ of
${\rm SL}(V)$-invariants is Frobenius split in this case also.
\elem \noindent {\it Proof.} Let~$\phi$ be a splitting of
$\spec(R)$. Continue to denote by~$\phi$ the restriction of~$\phi$
to the open part $\spec(R[1/f])$. We can `lift'~$\phi$ to the
$\mathbb{G}_m$-bundle $\spec(R[1/f][u,u^{-1}])$ (over
$\spec(R[1/f])$) as follows: define $\phi(a+\sum b_i u^i+\sum c_j u^{-j}) :=\phi(a)+\sum \phi(b_i)u^{i/p}+\sum \phi(c_j) u^{-j/p}$,  where the summations
are over positive integers and $u^{i/p}$ (respectively $u^{-j/p}$) is
interpreted to be~$0$ unless~$i$ (respectively $-j$) 
is an integral multiple of~$p$.
Observe that $R[1/f][u,u^{-1}]=S[1/f]$,  so we have a splitting of
$\spec(S[1/f])$ which we still denote~$\phi$.    We claim
that~$\phi$ maps~$S$ to~$S$ and hence extends to a splitting
of~$\spec(S)$.   Indeed,  a general element~$s$ of~$S$ is of the
form~$a+\sum b_i u^i+\sum c_j \xi^j$ with $a$,~$b_i$, and~$c_j$
in~$R$,   so that $\phi(s)=\phi(a+\sum b_iu^i+\sum c_j f^ju^{-j})
=\phi(a)+\sum \phi(b_i)u^{i/p}+\sum \phi(c_jf^j)u^{-j/p}$.
Rewriting $\phi(c_jf^j)u^{-j/p}$ as
$\phi(c_j)f^{j/p}u^{-j/p}=\phi(c_j)\xi^{j/p}$, we see that
$\phi(s)$ belongs to~$S$. \hfill $\Box$

\brem In the case when one of $\{m,q\}$ being $<n$, and the other
$\ge n$, we expect the ring of invariants to be Frobenius split
though at the moment, we do not have a proof of this assertion!
\erem


\noindent
{\bf Acknowledgments:}  Part of this work was carried out at
Abdus Salam International
Centre for Theoretical Physics, Trieste, during the visit of all the
three authors in April 2006. The authors gratefully acknowledge the financial
support and the hospitality of the ICTP.

\end{document}